\documentclass[11pt]{article}
\usepackage{amsfonts}
\usepackage{amsmath}
\usepackage{tabularx}
\usepackage{amssymb}
\usepackage{hyperref}
\usepackage{bbm}
\usepackage{microtype}
\newtheorem{thm}{Theorem}[section]

\newtheorem{defn}[thm]{Definition}
\newenvironment{defn-new}{\begin{defn} \em}{\end{defn}}
\newtheorem{rem}[thm]{Remark}
\newenvironment{rem-new}{\begin{rem} \em}{\end{rem}}
\newtheorem{ex}[thm]{Example}
\newenvironment{ex-new}{\begin{ex} \em}{\end{ex}}

\newenvironment{notation-new}{\begin{rem} \em}{\end{rem}}

\newenvironment{agr-new}{\begin{rem} \em}{\end{rem}}

\makeatletter \@addtoreset{equation}{section} \makeatother

\makeatletter \@addtoreset{figure}{section} \makeatother
\begin{document}
\begin{center}
{\textbf {\Large Explicit Constructions  of Astheno-K\"ahler Manifolds}}\\[2pt]
{\bf Punam Gupta}\footnote{
School of Mathematics, Devi Ahilya Vishwavidyalaya, Indore-452 001, M.P. 
INDIA\newline
Email: punam2101@gmail.com} and {\bf Nidhi Yadav}\footnote{
Department of Mathematics \& Statistics, Dr. Harisingh Gour Vishwavidyalaya,
Sagar-470 003, M.P. INDIA \newline
Email: nidhiyadav.bina@gmail.com} 
\end{center}

\noindent
{\bf Abstract:}
We investigate the conditions
under which astheno-K\"ahler structures can be identified on the product of two compact trans-Sasakian manifolds of dimensions greater than $2$.


\section{Introduction}

Over the past three decades, astheno-K\"ahler manifolds have garnered considerable attention from geometers globally, demonstrating their potential significance in both differential and algebraic frameworks. The phrase astheno-K\"ahler manifold, derived from the Greek word for weak, was introduced in 1993 by Jost and Yau in their influential paper \cite{jost,jost2}.

Matsuo and Takahashi \cite{km} demonstrated a pivotal finding, establishing that every compact balanced astheno-K\"ahler manifold is indeed K\"ahler. Notably, this implies that compact Hermitian-flat astheno-K\"ahler manifolds are also K\"ahler. Their research further revealed several important results:\begin{itemize}\item  any product manifold composed of curves and surfaces qualifies as astheno-Kähler.\item  for two compact Sasakian manifolds \( (M_i, g_i) \) in three dimensions, the product manifold \( M = M_1 \times M_2 \) equipped with the specified Hermitian structure and K\"ahler form is astheno-K\"ahler. 
\item if \( (M_1, g_1) \) is a three-dimensional compact Sasakian manifold and \( (M_2, g_2) \) is a compact cosymplectic manifold of dimension at least three, then the product manifold \( M = M_1 \times M_2 \) with the given Hermitian structure is also astheno-Kähler. 
\end{itemize}

In a separate study, Matsuo \cite{kmm} identified the existence of astheno-K\"ahler metrics on Calabi-Eckmann manifolds, which are principal \( T^2 \)-bundles over \({ \mathbbm CP}^n \times {\mathbbm CP}^m \), thereby facilitating the construction of astheno-K\"ahler structures on torus bundles over K\"ahler manifolds. 

Building on these findings, this work aims to comprehensively explore the conditions under which astheno-K\"ahler structures can be identified on various manifolds.

\section{Preliminaries}

 An almost complex manifold $(M,J)$ \cite{ballmann} together with a compatible
Riemannian metric $g$, that is, $g(JX, JY ) = g( X, Y)$
for all vector fields $X$, $Y$ on $M$, is called an almost Hermitian manifold and metric $g$ is called an almost Hermitian metric. A complex manifold $M$ \cite{ballmann} together with a compatible
Riemannian metric $g$ is called a Hermitian manifold and metric $g$ is called a Hermitian metric. 
The alternating $2$-form
 $$\Omega(X, Y ) = g(JX, Y)$$
is called the associated K\"{a}hler form. We can retrieve $g$ from $\Omega$
$$g( X, Y ) = \Omega(X, JY ).$$
\noindent If $\Omega $ is closed, then $(M,g)$ is known as a K\"{a}hler
manifold\cite{ballmann}  and $g$ is a K\"{a}hler metric.

%A quasi-K\"ahler structure is an almost Hermitian structure whose K\"ahler form $\Omega$ satisfies $$(d\Omega)^{(1,2)} = \bar{\partial} \Omega =0.$$
There are certain complex manifolds that do not support K\"ahler metrics, leading to the natural emergence of the question regarding suitable generalisations. While a universal type of Hermitian non-K\"ahler metrics has yet to be discovered, various classes associated with different geometric or physical applications have been proposed and examined. 

In 1993, Jost and Yau first coined the term astheno-K\"ahler manifold  in the remarkable paper \cite{jost,jost2}.

\begin{defn-new} {\rm\cite{jost,jost2}} A Hermitian manifold $(M,J,g)$ of complex dimension $m$ is called an astheno-K\"ahler manifold if it carries a fundamental $2$-form (K\"ahler form) $\Omega$ satisfying

 $$\partial \bar{\partial} \Omega^{m-2}=0,$$ where $\partial $
 and $\bar{\partial}$ are complex exterior differentials and $\Omega^{m-2}=\underset{m-2\quad \mathrm{times}}{\Omega \wedge \ldots\wedge \Omega}$.
\end{defn-new}
The above condition is automatically satisfied for $m=2$. Thus, any Hermitian metric on a complex surface is astheno-K\"ahler.

\begin{defn-new}
{\rm \cite{Streets}} Let  $(M,J,g)$  be a Hermitian manifold of complex dimension $m$. The Hermitian structure $(J,g)$ is called strong K\"ahler with torsion (strong KT or SKT)  and $g$ is called strong KT metric or SKT metric if $ \partial\Bar{\partial}\Omega = 0.$
\end{defn-new}

 \begin{defn-new}
{\rm \cite{gaud}} Let  $(M,J,g)$  be a Hermitian manifold of complex dimension $m$. The Hermitian structure $(J,g)$ is called standard or Gauduchon and $g$ is called Gauduchon metric if $ \partial\Bar{\partial}\Omega^{m-1} = 0$. Equivalently, if the Lee form of $\Omega^{m-1}$ is co-closed. 
 \end{defn-new}
\begin{rem-new}
    
Fino et al. {\rm \cite{fino}} mentioned that for $m = 4$, a Hermitian metric that is simultaneously strong KT and astheno-K\"ahler metric, it must also be Gauduchon. For $m=3$, astheno-K\"ahler structure means a strong KT metric.
\end{rem-new}
An odd-dimensional Riemannian manifold $\left(M^{2 n+1}, g\right)$ is said to be an almost contact metric manifold if there exists on $M$ a $(1,1)$ tensor field $\phi$, a vector field $\xi$ (called the structure vector field) and a 1-form $\eta$ such that

$$
\eta(\xi)=1, \quad \phi^2(X)=-X+\eta(X) \xi, \quad  \quad g(\phi X, \phi Y)=g(X, Y)-\eta(X)\eta(Y)
$$
for any vector fields $X, Y$ on $M$.\\

In particular, in an almost contact metric manifold, we also have $\phi \xi=0$ and $\eta \circ \phi=0$.
Such a manifold is said to be a contact metric manifold if $d \eta=\Phi$, where $\Phi(X, Y)=g(X, \phi Y)$ is called the fundamental $2 $-form of $M$. Such an almost contact metric structure is called a contact metric structure. Moreover, if a contact metric structure is normal, then it is called a Sasakian structure. On the other hand, if $d\Phi = 0$ and $d\eta = 0$, then $M$ is said to have an almost cosymplectic structure. In addition, if an almost cosymplectic structure is normal, then it is called a cosymplectic structure. If $N$ is a compact K\"ahler manifold, then $N \times S^1$ is the trivial example of compact cosymplectic manifolds. If a contact metric structure$(\phi, \xi, \eta)$ is normal and $d\eta = 0$, $d\Phi = 2\Phi \wedge \eta$, then $M$ is said to be Kenmotsu Manifold.

\noindent
These manifolds can be characterised through their Levi-Civita connection by requiring \\
(1): \text{Sasakian} $(\nabla_X \phi)Y = g(X,Y)\xi -\eta(Y)X,$ \\
(2): \text{Cosymplectic } $(\nabla_X \phi)Y = 0, $\\
(3): \text{Kenmotsu } $(\nabla_X \phi)Y = g(\phi X, Y)\xi -\eta(Y)\phi X.$

\begin{defn-new}{\rm \cite{janssens}} An almost $\alpha$-Sasakian manifold $M$ is an almost co-Hermitian manifold (or almost contact metric manifold) such that $\Phi=\left(1/\alpha\right) d \eta$, $\alpha \in {\mathbbm R}_0(\alpha \in{\mathbbm R},\alpha \neq 0 )$.$M$ is an $\alpha$-Sasakian manifold if the structure is normal.
\end{defn-new}

\begin{thm}{\rm \cite{janssens}}
{\rm(i)} An almost co-Hermitian manifold (or almost contact metric manifold) $M$ is $\alpha$-Sasakian if and only if 

$$
\left(\nabla_X \phi\right) Y=\alpha\{g(X, Y) \xi-\eta(Y) X\}
$$
for all $X, Y \in {{ \mathfrak X}}(M)$. 
\\{\rm(ii)} If $M$ is $\alpha$-Sasakian, then $\xi$ is a Killing vector field and

$$
\nabla_X \xi=-\alpha \phi X
$$
for all $X \in {\mathfrak{X}}(M)$.
\end{thm}
\begin{defn-new}{\rm \cite{janssens}} An almost $\beta$-Kenmotsu manifold is an almost co-Hermitian manifold (or almost contact metric manifold) such that

$$
\begin{gathered}
d \eta=0,\\
(d \Phi)(X, Y, Z)=\frac{2}{3} \beta \underset{X, Y, Z}{\mathfrak{S}}\{\eta(X) \Phi(Y, Z)\}
\end{gathered}
$$
for $\beta \in {\mathbbm R}_0$ and all $X, Y, Z \in {\mathfrak{X}}(M)$, where $\mathfrak{S}$ denotes the cyclic sum. $M$ is an $\beta$-Kenmotsu manifold if the structure is normal.
\end{defn-new}

\begin{thm}{\rm \cite{janssens}}
{\rm(i)} An almost co-Hermitıan manifold(or almost contact metric manifold) is an $\beta$-Kenmotsu manifold if and only if

$$
\left(\nabla_X \phi\right) Y=\beta\{g(\phi X, Y) \xi-\eta(Y) \phi X\},
$$
for all $X, Y \in {\mathfrak{X}}(M)$.
\end{thm}

\begin{defn-new}
{\rm \cite{oubina1}} A structure $(\phi, \xi, \eta, g)$ is a trans-Sasakian structure if and only if it is normal and

$$
d \eta=\alpha \Phi, \quad d \Phi=2 \beta \eta \wedge \Phi,
$$
where $\alpha=\displaystyle\frac{1}{2 n} \delta \phi(\xi)$, $\beta=\displaystyle\frac{1}{2 n}$ div $\xi$ and $\delta$ is the codifferential of $g$ and manifold is said to be trans-sasakian of type $(\alpha,\beta)$.\\
A trans-Sasakian structure of type $(\alpha,\beta)$ maybe expressed as an almost contact metric structure satisfying

$$(\nabla_{X} \phi) Y=\alpha(g(X, Y) \xi-\eta(Y) X)+\beta(g(\phi X, Y) \xi-\eta(Y) \phi X ).$$
\end{defn-new}
\noindent
From this formula, one easily obtains

\begin{equation}
\nabla_X \xi  =-\alpha \phi X-\beta \phi^2 X, \quad
\left(\nabla_X \eta\right) Y  =\alpha g(X, \phi Y)+\beta g(\phi X, \phi Y).
\end{equation}

\noindent
It is clear that a trans-Sasakian of type $(1,0)$ is a Sasakian manifold, a trans-Sasakian of type $(0,1)$ is a Kenmotsu manifold and a trans-Sasakian of type $(0,0)$ is a cosymplectic manifold.A $(\alpha,\beta)$ trans-sasakian manifold is an $\alpha$-Sasakian when $\alpha \in {\mathbbm R},\alpha \neq 0$ and $\beta= 0$ and is $\beta$-Kenmotsu when $\beta \in {\mathbbm R},\beta \neq 0$ and $\alpha= 0$.

\vspace{.2cm}
\noindent
 Now, consider the product manifold $M=M_1 \times M_2$, where $M_i(i=1,2)$ be a $\left(2 m_i+1\right)$-dimensional compact normal almost contact metric manifold with the structure tensor fields $\left(\phi_i, \xi_i, \eta_i\right)$,  we consider an almost complex structure $J$ (A.Morimoto's complex structure){\rm \cite{morimoto}} defined by
\begin{equation}\label{e3}
J=\phi_1-\eta_2 \otimes \xi_1+\phi_2+\eta_1 \otimes \xi_2 
\end{equation}
 Here $M$ endowed with $J$ is a compact complex manifold of complex dimension $m=m_1+m_2+1$.
\noindent
Moreover, if $g_i$ is the compatible Riemannian metric on $M_i$ for each $i=1,2$, then the Riemannian product metric $g=g_1+g_2$ on $M$ is compatible with $J$, that is, $g$ is a Hermitian metric on $M$. 
\noindent
Then its K\"ahler form $\Omega$  is given by
\begin{equation}\label{e2}
\Omega=\Phi_1+\Phi_2-2 \eta_1 \wedge \eta_2,  
\end{equation}
where $\Phi_i$ denotes the fundamental 2-form on $M_i$ for each $i=1,2$.

Marrero \cite{marrero} proved the following two results for trans-Sasakian manifolds:
\begin{thm}
  A $(\alpha,\beta)$ trans-sasakian manifold of dim $m=2n+1\geq 5$ is either $\alpha$-Sasakian, $\beta$-Kenmotsu or cosymplectic.
  \label{t1}
\end{thm}

\section{Construction of astheno-K\"ahler manifolds}
 Let $(M_i,\phi_i, \xi_i, \eta_i, g_i)$, $i=1,2$ be two compact trans-Sasakian manifolds of type and dimension $(\alpha_i,\beta_i)$ and $2m_i+1$, respectively, where $\alpha_i,\beta_i $ are functions on $M_i$ and  $\alpha_i=\displaystyle\frac{1}{2 n_i} \delta_i \phi_i(\xi_i)$, $\beta_i=\displaystyle\frac{1}{2 n_i} \operatorname{div} \xi_i$, $ \delta_i$ is the codifferential of $g_i$.
 
 \vspace{.2 cm}
 \noindent
 Let $M=M_1 \times M_2$ be a product manifold with complex dimension $ m=m_1+m_2+1$.
 
Consider the K\"ahler form $\Omega$ 
\begin{eqnarray}
\Omega &=&\Phi_1+\Phi_2-2 \eta_1 \wedge \eta_2. 
\label{p1}
\end{eqnarray}
By  (\ref{p1}), we get 
\begin{eqnarray}d \Omega &=& d \Phi_1+d \Phi_2-2\left(d \eta_1 \wedge \eta_2+\eta_1 \wedge d \eta_2\right), \notag \\
d \Omega&=&2 \beta_1 \eta_1 \wedge \Phi_1+2 \beta_2 \eta_2\wedge \Phi_2-2\left(\alpha_1 \Phi_1 \wedge \eta_2+\eta_1 \wedge \alpha_2 \Phi_2\right),\label{Eq:result} 
\end{eqnarray}

and 
\begin{eqnarray}d^c \Omega&=& Jd \Omega=2\left(J \beta_1 \eta_2 \wedge \Phi_1-J \beta_2 \eta_1 \wedge \Phi_2+J \alpha_1 \Phi_1 \wedge \eta_1\right. -\left.\eta_2 \wedge J \alpha_2 \Phi_2\right), \notag\\
&=&2\left(\beta_1 \eta_2 \wedge \Phi_1-\beta_2 \eta_1 \wedge \Phi_2+\alpha_1 \Phi_1 \wedge \eta_1-\alpha_2 \eta_2 \wedge \Phi_2\right). 
\end{eqnarray}
Then 
\begin{align}
 dd^c \Omega =&2\{\beta_1 d \eta_2 \wedge \Phi_1+\beta_1 \eta_2 \wedge d \Phi_1-\beta_2 d \eta_1 \wedge \Phi_2-\beta_2 \eta_1 \wedge d \Phi_2+\alpha_1 d \Phi_1 \wedge \eta_1 \notag \\
 &+\alpha_1 \Phi_1 \wedge d \eta_1-\alpha_2 d \eta_1 \wedge \Phi_2-\alpha_2 \eta_2 \wedge d \Phi_2\} \notag \\
 =&2\{\beta_1 \alpha_2 \Phi_2 \wedge \Phi_1+\beta_1 \eta_2 \wedge 2 \beta_1 \eta_1 \wedge \Phi_1-\beta_2 \alpha_1 \Phi_1 \wedge \Phi_2 -\beta_2 \eta_1 \wedge 2 \beta_2 \eta_2 \wedge \Phi_2 \notag \\
&+2 \alpha_1 \beta_1 \eta_1 \wedge \Phi_1\wedge \eta_1+\alpha_1 \Phi_1 \wedge \alpha_1 \Phi_1 -\alpha_2 \alpha_2 \Phi_2 \wedge \Phi_2-\alpha_2 \eta_2 \wedge 2 \beta_2 \eta_2 \wedge \Phi_2\} \notag \\
    = &2\{\beta_1 \alpha_2 \Phi_2 \wedge \Phi_1+2 \beta_1^2 \eta_2 \wedge \eta_1 \wedge \Phi_1-\beta_2 \alpha_1 \Phi_1 \wedge \Phi_2-2 \beta_2^2 \eta_1 \wedge \eta_2 \wedge \Phi_2  +\alpha_1^2 \Phi_1^2-\alpha_2^2 \Phi_2^2 \}.
    \label{c1}
\end{align}
Now
\begin{align}
 d \Omega \wedge d^c \Omega=&4\left(\beta_1 \eta_1 \wedge \Phi_1+\beta_2 \eta_2 \wedge \Phi_2-\alpha_1 \Phi_1 \wedge \eta_2-\eta_1 \alpha_2 \wedge \Phi_2\right) \wedge \notag \\
& \left(\beta_1 \eta_2 \wedge \Phi_1-\beta_2 \eta_1 \wedge \Phi_2+\alpha_1 \Phi_1 \wedge \eta_1-\alpha_2 \eta_2 \wedge \Phi_2\right) \\
 =&4[\beta_1^2 \eta_1 \wedge \Phi_1 \wedge \eta_2 \wedge \Phi_1-\beta_1 \beta_2 \eta_1 \wedge \Phi_1 \wedge \eta_1 \wedge \Phi_2 +\beta_1 \alpha_1 \eta_1 \wedge \Phi_1 \wedge \Phi_1 \wedge \eta_1 \notag \\
&-\beta_1 \alpha_2 \eta_1 \wedge \Phi_1 \wedge \eta_2 \Phi_2  +\beta_2 \beta_1 \eta_2 \wedge \eta_2 \wedge \Phi_2 \wedge \Phi_1-\beta_2^2 \eta_2 \wedge \Phi_2 \wedge \eta_1 \wedge \Phi_2 \notag  \\
& +\beta_2 \alpha_1 \eta_2 \wedge \Phi_2 \wedge \Phi_1 \wedge \eta_1-\beta_2 \alpha_2 \eta_2 \wedge \Phi_2 \wedge \eta_2 \wedge \Phi_2  -\alpha_1 \beta_1 \Phi_1 \wedge \eta_2 \wedge \eta_2 \wedge \Phi_1 \notag \\
&+\alpha_1 \beta_2 \Phi_1 \wedge \eta_2 \wedge \eta_1 \wedge \Phi_2  -\alpha_1^2 \Phi_1 \wedge \eta_2 \wedge \Phi_1 \wedge \eta_1+\alpha_1 \alpha_2 \Phi_1 \wedge \eta_2 \wedge \eta_2 \wedge \Phi_2 \notag \\
& -\beta_1 \alpha_2 \eta_1 \wedge \Phi_2 \wedge \eta_2 \wedge \Phi_1+\alpha_2 \beta_2 \eta_1 \wedge \Phi_2 \wedge \eta_1 \wedge \Phi_2 -\alpha_2 \alpha_1 \eta_1 \wedge \Phi_2 \wedge \Phi_1 \wedge \eta_1 \notag \\
&+\alpha_2^2 \eta_1 \wedge \Phi_2 \wedge \eta_2 \wedge \Phi_2] \notag \\
=&4[\beta_1^2 \eta_1 \wedge \eta_2 \wedge \Phi_1^2-\beta_1 \alpha_2 \eta_1 \wedge \eta_2 \wedge \Phi_1 \wedge \Phi_2 -\beta_2^2 \eta_2 \wedge \eta_1 \wedge \Phi_2^2+\beta_2 \alpha_1 \eta_2 \wedge \eta_1\wedge\Phi_2 \wedge \Phi_1 \notag \\
& +\alpha_1 \beta_2 \eta_2 \wedge \eta_1 \wedge \Phi_1 \wedge \Phi_2-\alpha_1^2 \eta_2 \wedge \eta_1 \wedge \Phi_1^2 -\beta_1 \alpha_2 \eta_1 \wedge \eta_2 \wedge \Phi_2 \wedge \Phi_1+\alpha_2^2 \eta_1\wedge \eta_2 \wedge \Phi_2^2].
\label{p6}
\end{align}
On the other hand, we have 
\begin{eqnarray}
 dd^c \Omega \wedge \Omega&=&2\left\{\beta_1 \alpha_2 \Phi_2 \wedge \Phi_1+2 \beta_1^2 \eta_2 \wedge \eta_1 \wedge \Phi_1-\beta_2 \alpha_1 \Phi_1 \wedge \Phi_2 -2 \beta_2^2 \eta_1 \wedge \eta_2 \wedge \Phi_2+\alpha_1^2 \Phi_1^2-\alpha_2^2 \Phi_2^2\right\} \notag \\
&& \wedge\left\{\Phi_1+\Phi_2-2 \eta_1 \wedge \eta_2\right\} \notag\\
&=&2[\beta_1 \alpha_2 \Phi_2 \wedge \Phi_1^2+\beta_1 \alpha_2 \Phi_2^2 \wedge\Phi_1-2 \beta_1 \alpha_2 \Phi_2 \wedge \Phi_1 \wedge \eta_1 \wedge \eta_2 +2 \beta_1^2 \eta_2 \wedge \eta_1 \wedge \Phi_1^2 \notag \\
&&+2 \beta_1^2 \eta_2 \wedge \eta_1 \wedge \Phi_1 \wedge \Phi_2 -4 \beta_1^2 \eta_2 \wedge \eta_1\wedge \Phi_1\wedge \eta_1 \wedge \eta_2  -\beta_2 \alpha_1 \Phi_1^2 \wedge \Phi_2-\beta_2 \alpha_1 \Phi_1 \wedge \Phi_2^2 \notag\\
&&+2 \beta_2 \alpha_1 \Phi_1 \wedge \Phi_2 \wedge \eta_1 \wedge \eta_2 -2 \beta_2^2 \eta_1 \wedge \eta_2 \wedge \Phi_2 \wedge \Phi_1-2 \beta_2^2 \wedge \eta_1 \wedge \eta_2 \wedge \Phi_2^{2}\notag  \\
&& +4 \beta_2^2 \eta_1 \wedge \eta_2 \wedge \Phi_2 \wedge \eta_1 \wedge \eta_2 +\alpha_1^2 \Phi_1^3 +\alpha_1^2 \Phi_1^2 \wedge \Phi_2-2 \alpha_1^2 \Phi_1^2 \wedge \eta_1 \wedge \eta_2 -\alpha_2^2 \Phi_2^2 \wedge \Phi_1\notag \\
&& -\alpha_2^2 \Phi_2^3+2 \alpha_2^2 \Phi_2^2 \wedge \eta_1 \wedge \eta_2].
\label{p7}
\end{eqnarray}
If complex dimension of manifold is $m$, then for astheno-K\"ahler we need to show

\begin{eqnarray}
dd^c \Omega^{m-2} = (m-2)[dd^c\Omega \wedge \Omega + (m-3)d\Omega \wedge d^c \Omega]\wedge \Omega^{m-4}=0.
\label{e1}
 \end{eqnarray}
 The possibilities of dimension $m$ of the product manifold $M$ (over $\mathbbm{C}$) are following:\\
\noindent
\textbf{Case 1:} For $m = 3$ (complex dimension), that is, we take the product of two $3$-trans-Sasakian manifolds, that is, the dimension of $M_1$ is $ 2m_1 + 1 = 3$ and dimension of $M_2$ is $2m_2 + 1=3$. To show a manifold of dimension $m=3$ is astheno-K\"ahler, we need to prove $dd^c \Omega=0$.
From (\ref{c1}),\\
$dd^c \Omega = 2 \left\{ \beta_1 \alpha_2 \Phi_2 \wedge \Phi_1 + 2 \beta_1^2 \eta_2 \wedge \eta_1 \wedge \Phi_1 -\beta_2 \alpha_1 \wedge\Phi_1 \wedge\Phi_2 -2 \beta_2^2 \eta_1 \wedge \eta_2 \wedge \Phi_2 + \alpha_1^2 \Phi_1^2 -\alpha_2^2 \Phi_2^2 \right\}$

\vspace{.5cm}
\noindent
\begin{tabular}{|p{.7 cm}|p{2.5cm}|p{2.5cm}|p{.2 cm}|p{.2 cm}|p{.2 cm}|p{.2 cm}|p{6.6cm}| }
 \hline
 \multicolumn{8}{|c|}{Table 1: Particular cases} \\
 \hline
 S.No. &$M_1$ & $M_2$& $\alpha_1$ & $\alpha_2$ & $\beta_1$&$\beta_2$& 
 $dd^c \Omega $ \\
 \hline
1&$\alpha_1$-Sasakian & $\alpha_2$-Sasakian  & 1 & 1 &0  & 0 & 0 \\
\hline
2&$\alpha_1$-Sasakian & $\beta_2$-Kenmotsu  & 1 &0  & 0 &1  &
$-2  \beta_2 \alpha_1\Phi_1 \wedge \Phi_2 $ 
$-4 \beta_2^2 \eta_1 \wedge \eta_2 \wedge \Phi_2 +2\alpha_1^2\Phi_1^2$ \\
\hline
3&$\alpha_1$-Sasakian & Cosymplectic   &1  & 0 & 0 &0  &0   \\
\hline
4&$\beta_1$-Kenmotsu & $\beta_2$-Kenmotsu  & 0 & 0 & 1 &1  & 
$2 \{ 2 \beta_1^2 \eta_2 \wedge \eta_1 \wedge \Phi_1 $ 
$-2 \beta_2^2 \eta_1 \wedge \eta_2 \wedge \Phi_2 $\}  \\
\hline
5&$\beta_1$-Kenmotsu & $\alpha_2$-Sasakian &0  & 1  & 1 & 0 & 
$2 \{ \beta_1 \alpha_2 \Phi_2 \wedge \Phi_1$ 
$ + 2 \beta_1^2 \eta_2 \wedge \eta_1 \wedge \Phi_1 \}$    \\
\hline
6&$\beta_1$-Kenmotsu & Cosymplectic & 0 & 0 & 1 & 0 & $2 \beta_1^2 \eta_2 \wedge \eta_1 \wedge \Phi_1$    \\
\hline
7&Cosymplectic & $\alpha_2$-Sasakian  & 0 & 1 & 0 & 0 & 0 \\
\hline
8&Cosymplectic &  $\beta_2$-Kenmotsu  & 0 & 0 & 0 & 1 & $-2 \beta_2^2 \eta_1 \wedge \eta_2 \wedge \Phi_2$   \\
\hline
9&Cosymplectic &   Cosymplectic  & 0 & 0 & 0 & 0 & $0$  \\
 \hline
 
\end{tabular}
\begin{rem-new}
    1. In Table 1, Cases 1, 3, 7 and 9 are astheno-K\"ahler.\\
    2. Specific instances of Case 1 and Cases 3 and 7 can be observed in Theorems 4.1 and 4.2, respectively, in the paper by Matsuo and Takahashi \cite{km}.

\end{rem-new}

\vspace{.5cm}
\noindent
\textbf{Case 2 }: For $m = 4$ (complex dimension), there are two possibilities.

\noindent (i) One of possibility is $m_1 = 1$, $m_2 = 2$, that is, $\dim M_1 = 3$, $\dim M_2 = 5$. From (\ref{e1}), we have
$dd^c \Omega^2 = 2 \left[ dd^c \Omega \wedge \Omega + d\Omega \wedge d^c \Omega \right]=0. $

\noindent Using (\ref{p6}) and (\ref{p7}), we have
\begin{align*}
dd^c \Omega^2 =& 2 \bigg[ 2 ( \beta_1 \alpha_2 \Phi_2 \wedge \Phi_1^2 + \beta_1 \alpha_2 \Phi_2^2 \wedge \Phi_1  -2 \beta_1 \alpha_2 \Phi_2 \wedge \Phi_1\wedge\eta_1 \wedge \eta_2 + 2 \beta_1^2 \eta_2 \wedge \eta_1 \wedge \Phi_1^2 \\
& + 2 \beta_1^2 \eta_2 \wedge \eta_1 \wedge \Phi_1\wedge \Phi_2 -\beta_2 \alpha_1 \wedge\Phi_1^2\wedge \Phi_2 -\beta_2 \alpha_1 \Phi_1\wedge \Phi_2^2   \\
&+ 2 \beta_2 \alpha_1\wedge\Phi_1\wedge\Phi_2 \wedge \eta_1 \wedge \eta_2   -2 \beta_2^2 \eta_1 \wedge \eta_2 \wedge\Phi_2\wedge \Phi_1 -2 \beta_2^2 \eta_1 \wedge \eta_2 \wedge \Phi_2^2 \\
& + \alpha_1^2\Phi_1^3 + \alpha_1^2 \Phi_1^2 \wedge \Phi_2 -2 \alpha_1^2 \Phi_1^2 \wedge \eta_1 \wedge \eta_2 -\alpha_2^2 \Phi_2^2 \wedge \Phi_1 -\alpha_2^2 \Phi_2^3 + 2 \alpha_2^2 \Phi_2^2 \wedge \eta_1 \wedge \eta_2 )\\
&+4(\beta_1^2 \eta_1 \wedge \eta_2 \wedge \Phi_1^2-\beta_1 \alpha_2 \eta_1 \wedge \eta_2 \wedge \Phi_1 \wedge \Phi_2 -\beta_2^2 \eta_2 \wedge \eta_1 \wedge \Phi_2^2+\beta_2 \alpha_1 \eta_2 \wedge \eta_1\wedge\Phi_2 \wedge \Phi_1 \\
& +\alpha_1 \beta_2 \eta_2 \wedge \eta_1 \wedge \Phi_1 \wedge \Phi_2-\alpha_1^2 \eta_2 \wedge \eta_1 \wedge \Phi_1^2 -\beta_1 \alpha_2 \eta_1 \wedge \eta_2 \wedge \Phi_2 \wedge \Phi_1+\alpha_2^2 \eta_1\wedge \eta_2 \wedge \Phi_2^2) \bigg]\\
=& 2 \Big[ 2 \beta_1 \alpha_2 \Phi_2\wedge\Phi_1^2 + 2 \beta_1 \alpha_2 \Phi_2^2\wedge \Phi_1 -12 \beta_1\alpha_2 \eta_1\wedge \eta_2 \wedge\Phi_2\wedge \Phi_1  + 4 \beta_1^2 \eta_1 \wedge\eta_2 \wedge\Phi_1 \wedge\Phi_2 \\
&-2 \beta_2 \alpha_1 \Phi_1^2\wedge \Phi_2 -2 \beta_2 \alpha_1 \Phi_1 \wedge\Phi_2^2  + 4 \alpha_1 \beta_2 \eta_2\wedge \eta_1\wedge \Phi_1 \wedge\Phi_2 -4 \beta_2^2 \eta_1\wedge \eta_2\wedge \Phi_2\wedge \Phi_1  \\
&+ 2 \alpha_1^2\Phi_1^3 + 2\alpha_1^2\Phi_1^2\wedge\Phi_2 -2 \alpha_2^2 \Phi_2^2 \wedge\Phi_1 -2 \alpha_2^2 \Phi_2^3 + 8 \alpha_2^2\eta_1\wedge \eta_2\wedge \Phi_2^2\Big]\quad .
\end{align*}

\noindent
\begin{tabular}{|p{.8 cm}|p{2.5cm}|p{2.5cm}|p{.2 cm}|p{.2 cm}|p{.2 cm}|p{.2 cm}|p{6.6cm}|  }
 \hline
 \multicolumn{8}{|c|}{Table 2: Particular cases} \\
 \hline
 S.No. &$M_1$ & $M_2$& $\alpha_1$ & $\alpha_2$ & $\beta_1$&$\beta_2$& $dd^c \Omega^2 $ \\
 \hline
1&$\alpha_1$-Sasakian & $\alpha_2$-Sasakian  & 1 & 1 &0  & 0 & $ 16 \alpha_2^2 \eta_1 \wedge \eta_2 \wedge \Phi_2^2 -4\alpha_2^2 \Phi_2^2 \wedge \Phi_1-4\alpha_2^2\Phi_2^2 $   \\
 \hline
2&$\alpha_1$-Sasakian & $\beta_2$-Kenmotsu  & 1 &0  & 0 &1  & $ 8\beta_2 \alpha_1 \Phi_1\wedge\Phi_2\wedge\eta_2\wedge\eta_1 -4\beta_2 \alpha_1 \Phi_1 \wedge \Phi_2^2-8 \beta_2^2 \eta_1 \wedge \eta_2 \wedge \Phi_2\wedge\Phi_1$ \\
 \hline
3&$\alpha_1$-Sasakian & Cosymplectic   &1  & 0 & 0 &0  & 0  \\
 \hline
4&$\beta_1$-Kenmotsu & $\beta_2$-Kenmotsu  & 0 & 0 & 1 &1& $  8 \beta_1^2 \eta_2 \wedge \eta_1 \wedge \Phi_1 \wedge \Phi_2 -8 \beta_2^2 \eta_1 \wedge \eta_2 \wedge \Phi_2\wedge\Phi_1$ \\
 \hline
5&$\beta_1$-Kenmotsu & $\alpha_2$-Sasakian &0  & 1  & 1 & 0 & $ -24 \beta_1 \alpha_2 \Phi_1 \wedge \Phi_2 \wedge \eta_1 \wedge \eta_2 + 16 \alpha_2^2 \eta_1 \wedge \eta_2 \wedge \Phi_2^2 +4 \beta_1 \alpha_2 \Phi_2^2 \wedge \Phi_1 + 8 \beta_1^2 \eta_2 \wedge \eta_1 \wedge \Phi_1\wedge\Phi_2 -4\alpha_2^2 \Phi_2^2 \wedge \Phi_1 +4\beta_1\alpha_2\Phi_2\wedge\Phi_1^2$ \\
\hline
6&$\beta_1$-Kenmotsu & Cosymplectic & 0 & 0 & 1 & 0 & $8 \beta_1^2 \eta_2 \wedge \eta_1 \wedge \Phi_1 \wedge \Phi_2$  \\
\hline
7&Cosymplectic & $\alpha_2$-Sasakian  & 0 & 1 & 0 & 0 & $16 \alpha_2^2 \eta_1 \wedge \eta_2 \wedge \Phi_2^2 -4\alpha_2^2 \Phi_2^2 \wedge \Phi_1$ \\
\hline
8&Cosymplectic &  $\beta_2$-Kenmotsu  & 0 & 0 & 0 & 1 & $ -8 \beta_2^2 \eta_1 \wedge \eta_2 \wedge \Phi_2 \wedge \Phi_1$ \\
 \hline
 9&Cosymplectic &  Cosymplectic  & 0 & 0 & 0 & 0 & 0 \\
 \hline
\end{tabular}
\begin{rem-new}
    1. In Table $2$, cases $3$ and $9$ are the astheno-K\"ahler manifold.\\
    2. A specific instance of Case 3 can be observed in Theorem 4.2 in the paper by Matsuo and Takahashi \cite{km}.
\end{rem-new}
\vspace{.3cm}
\noindent
(ii) Another possibility is $\dim M_1 = 5$, $\dim M_2 = 3$,

Since we have,
\begin{align*}
    dd^c \Omega^2 =& 2 \Big[ 2 \beta_1 \alpha_2 \Phi_2\wedge\Phi_1^2 + 2 \beta_1 \alpha_2 \Phi_2^2\wedge \Phi_1 -12 \beta_1\alpha_2 \eta_1\wedge \eta_2 \wedge\Phi_2\wedge \Phi_1  + 4 \beta_1^2 \eta_1 \wedge\eta_2 \wedge\Phi_1 \wedge\Phi_2 \\
&-2 \beta_2 \alpha_1 \Phi_1^2\wedge \Phi_2 -2 \beta_2 \alpha_1 \Phi_1 \wedge\Phi_2^2  + 4 \alpha_1 \beta_2 \eta_2\wedge \eta_1\wedge \Phi_1 \wedge\Phi_2 -4 \beta_2^2 \eta_1\wedge \eta_2\wedge \Phi_2\wedge \Phi_1  \\
& + 2 \alpha_1^2\Phi_1^3 + 2\alpha_1^2\Phi_1^2\wedge\Phi_2 -2 \alpha_2^2 \Phi_2^2 \wedge\Phi_1 -2 \alpha_2^2 \Phi_2^3 + 8 \alpha_2^2\eta_1\wedge \eta_2\wedge \Phi_2^2 \Big]
\end{align*}

\vspace{.5cm}
\noindent
\begin{tabular}{|p{.8 cm}|p{2.5cm}|p{2.5cm}|p{.2 cm}|p{.2 cm}|p{.2 cm}|p{.2 cm}|p{6.8 cm}|  }
 \hline
 \multicolumn{8}{|c|}{Table 3: Particular cases} \\
 \hline
 S.No. &$M_1$ & $M_2$& $\alpha_1$ & $\alpha_2$ & $\beta_1$&$\beta_2$& $dd^c \Omega^2 $ \\
 \hline
1&$\alpha_1-$ Sasakian & $\alpha_2$-Sasakian  & 1 & 1 &0  & 0 & $4\alpha_1^2\Phi_1^2\wedge\Phi_2$ \\
\hline
2&$\alpha_1$-Sasakian & $\beta_2$-Kenmotsu  & 1 &0  & 0 &1  & $8\beta_2 \alpha_1 \Phi_1\wedge\Phi_2\wedge\eta_2\wedge\eta_1 -4\beta_2 \alpha_1 \Phi_1 \wedge \Phi_2^2-8 \beta_2^2 \eta_1 \wedge \eta_2 \wedge \Phi_2\wedge\Phi_1 -4\beta_2\alpha_1\Phi_1^2\wedge\Phi_2 +4\alpha_1^2\Phi_1^2\wedge\Phi_2$ \\
\hline
3&$\alpha_1$-Sasakian & Cosymplectic   &1  & 0 & 0 &0  & $4\alpha_1^2\Phi_1^2\wedge\Phi_2$  \\
\hline
4&$\beta_1$-Kenmotsu & $\beta_2$-Kenmotsu  & 0 & 0 & 1 &1 & $ -8 \beta_2^2 \eta_1 \wedge \eta_2 \wedge \Phi_2\wedge\Phi_1 + 8 \beta_1^2 \eta_1 \wedge \eta_2 \wedge \Phi_1 \wedge \Phi_2 $   \\
\hline
5&$\beta_1$-Kenmotsu & $\alpha_2$-Sasakian &0  & 1  & 1 & 0& $  4\beta_1 \alpha_2 \Phi_2\wedge\Phi_1^2 -24\beta_1 \alpha_2 \Phi_1 \wedge \Phi_2\wedge\eta_1 \wedge \eta_2 + 8 \beta_1^2 \eta_1 \wedge \eta_2 \wedge \Phi_1 \wedge \Phi_2   $   \\
\hline
6&$\beta_1$-Kenmotsu & Cosymplectic & 0 & 0 & 1 & 0 &$  8\beta_1^2 \eta_1 \wedge \eta_2 \wedge \Phi_1\wedge\Phi_2 $  \\
\hline
7&Cosymplectic & $\alpha_2$-Sasakian  & 0 & 1 & 0 & 0 & $ 0  $ \\
\hline
8&Cosymplectic &  $\beta_2$-Kenmotsu  & 0 & 0 & 0 & 1 & $ -8\beta_2^2 \eta_1 \wedge \eta_2 \wedge \Phi_2 \wedge \Phi_1 $ \\
 \hline
 9 &Cosymplectic &  Cosymplectic  & 0 & 0 & 0 & 0 & 0 \\
 \hline
\end{tabular}

\begin{rem-new}
   1. In Table 3, Cases 7 and 9 are the astheno-K\"ahler manifold.\\
   2. A specific instance of Case 7 can be observed in Theorem 4.2 in the paper by Matsuo and Takahashi \cite{km}.
\end{rem-new}
\begin{rem-new}
    In case 4 of Tables 1,2 and 3, if $\beta_1 = \pm \beta_2$, it will become the astheno-K\"ahler manifold.
\end{rem-new}
\textbf{Case 3}: For $\dim_{\mathbbm{C}} M = m = 5$, that is, $\dim_{\mathbbm{R}} M = 10$. We have the following possibilities :
\begin{itemize}
    \item[(i)] $\dim_{\mathbbm{R}} M_1 = 3$ and $\dim_{\mathbbm{R}} M_2 = 7$
\item[(ii)] $\dim_{\mathbbm{R}} M_1 = 5$ and $\dim_{\mathbbm{R}} M_2 = 5$
\item[ (iii)] $\dim_{\mathbbm{R}} M_1 = 7$ and $\dim_{\mathbbm{R}} M_2 = 3$
\end{itemize}

\noindent
For astheno-K\"ahler manifold with $m = 5$ , we need to show 
$$dd^c \Omega^3 =3(dd^c \Omega \wedge \Omega + 2d \Omega \wedge d^c \Omega) \wedge \Omega = 0.$$ 
By (\ref{p6}) and (\ref{p7}),we have
\begin{eqnarray*}
 dd^c \Omega^3 &=&3\Big[2(\beta_1 \alpha_2 \Phi_2 \wedge \Phi_1^2+\beta_1 \alpha_2 \Phi_2^2 \wedge\Phi_1-2 \beta_1 \alpha_2 \Phi_2 \wedge \Phi_1 \wedge \eta_1 \wedge \eta_2 +2 \beta_1^2 \eta_2 \wedge \eta_1 \wedge \Phi_1^2  \\
&&+2 \beta_1^2 \eta_2 \wedge \eta_1 \wedge \Phi_1 \wedge \Phi_2  -\beta_2 \alpha_1 \Phi_1^2 \wedge \Phi_2-\beta_2 \alpha_1 \Phi_1 \wedge \Phi_2^2 \\
&&+2 \beta_2 \alpha_1 \Phi_1 \wedge \Phi_2 \wedge \eta_1 \wedge \eta_2 -2 \beta_2^2 \eta_1 \wedge \eta_2 \wedge \Phi_2 \wedge \Phi_1-2 \beta_2^2 \wedge \eta_1 \wedge \eta_2 \wedge \Phi_2^{2}  \\
&& +\alpha_1^2 \Phi_1^3 +\alpha_1^2 \Phi_1^2 \wedge \Phi_2-2 \alpha_1^2 \Phi_1^2 \wedge \eta_1 \wedge \eta_2 -\alpha_2^2 \Phi_2^2 \wedge \Phi_1 \\
&& -\alpha_2^2 \Phi_2^3+2 \alpha_2^2 \Phi_2^2 \wedge \eta_1 \wedge \eta_2) + 2\times4(\beta_1^2 \eta_1 \wedge \eta_2 \wedge \Phi_1^2-\beta_1 \alpha_2 \eta_1 \wedge \eta_2 \wedge \Phi_1 \wedge \Phi_2 \\
&&-\beta_2^2 \eta_2 \wedge \eta_1 \wedge \Phi_2^2 +\beta_2 \alpha_1 \eta_2 \wedge \eta_1\wedge\Phi_2 \wedge \Phi_1  +\alpha_1 \beta_2 \eta_2 \wedge \eta_1 \wedge \Phi_1 \wedge \Phi_2-\alpha_1^2 \eta_2 \wedge \eta_1 \wedge \Phi_1^2 \\
&&-\beta_1 \alpha_2 \eta_1 \wedge \eta_2 \wedge \Phi_2 \wedge \Phi_1+\alpha_2^2 \eta_1\wedge \eta_2 \wedge \Phi_2^2) \Big]\wedge\Big[\Phi_1+\Phi_2-2 \eta_1 \wedge \eta_2\Big] 
\end{eqnarray*}
\begin{eqnarray*}
    &=&3\Big[2\beta_1\alpha_2\Phi_2\wedge\Phi_1^3+ 4\beta_1\alpha_2\Phi_2^2\wedge\Phi_1^2- 24\beta_1\alpha_2\Phi_2\wedge\Phi_1^2\wedge\eta_1\wedge\eta_2+ 4\beta_1^2\eta_1\wedge\eta_2\wedge\Phi_1^3\\
   &&- 2\beta_2\alpha_1\Phi_1^3\wedge\Phi_2- 4\beta_2\alpha_1\Phi_1^2\Phi_2^2- 4\beta_2^2\eta_1\wedge\eta_2\wedge\Phi_2\wedge\Phi_1^2+ 2\alpha_1^2\Phi_1^4+ 4\alpha_1^2\Phi_1^4\\
   &&+ 4\alpha_1^2\Phi_1^3\wedge\Phi_2- 4\alpha_2^2\Phi_2^3\wedge\Phi_1- 2\alpha_2^2\Phi_2^2\wedge\Phi_1^2+ 16\alpha_2^2\eta_1\wedge\eta_2\wedge\Phi_2^2\wedge\Phi_1+ 2\beta_1\alpha_2\Phi_2^3\wedge\Phi_1\\
   &&- 24\beta_1\alpha_2\Phi_2^2\wedge\Phi_1\wedge\eta_1\wedge\eta_2+ 4\beta_1^2\eta_2\wedge\eta_1\wedge\Phi_1\wedge\Phi_2^2- 2\beta_2\alpha_1\Phi_1\wedge\Phi_2^3\\
   &&+ 8\beta_2\alpha_1\Phi_1\wedge\Phi_2^2\wedge\eta_2\wedge\eta_1+ 4\beta_2^2\eta_1\wedge\eta_2\wedge\Phi_2^3 + 2\alpha_1^2\Phi_1^2\wedge\Phi_2^2- 2\alpha_2^2\Phi_2^4+ 16\alpha_2^2\Phi_2^3\wedge\eta_1\wedge\eta_2 \Big] \quad.
\end{eqnarray*}

(i) $\dim_{\mathbbm{R}} M_1 = 3$ and $\dim_{\mathbbm{R}} M_2 = 7$

\vspace{.5cm}
\noindent
\begin{tabular}{|p{.8 cm}|p{2.5cm}|p{2.5cm}|p{.2 cm}|p{.2 cm}|p{.2 cm}|p{.2 cm}|p{6.8 cm}|  }
 \hline
 \multicolumn{8}{|c|}{Table 4: Particular cases} \\
 \hline
 S.No. &$M_1$ & $M_2$& $\alpha_1$ & $\alpha_2$ & $\beta_1$&$\beta_2$& $dd^c \Omega^3 $ \\
 \hline
1&$\alpha_1-$ Sasakian & $\alpha_2$-Sasakian  & 1 & 1 &0  & 0 &$3(-4\alpha_2^2\Phi_2^3\wedge\Phi_1+ 16\alpha_2^2\eta_1\wedge\eta_2\wedge\Phi_2^2\wedge\Phi_1+ 16\alpha_2^2\Phi_2^3\wedge\eta_1\wedge\eta_2)$ \\
\hline
2&$\alpha_1$-Sasakian & $\beta_2$-Kenmotsu  & 1 &0  & 0 &1 & $3(-2\beta_2\alpha_1\Phi_1\wedge\Phi_2^3+ 8\beta_2\alpha_1\Phi_1\wedge\Phi_2^2\wedge\eta_2\wedge\eta_1+ 4\beta_2^2\eta_1\wedge\eta_2\wedge\Phi_2^3)$ \\
\hline
3&$\alpha_1$-Sasakian & Cosymplectic   &1  & 0 & 0 &0  & $0$ \\
\hline
4&$\beta_1$-Kenmotsu & $\beta_2$-Kenmotsu  & 0 & 0 & 1 &1 & $3(4\beta_1^2\eta_1\wedge\eta_2\wedge\Phi_1^3- 4\beta_2^2\eta_1\wedge\eta_2\wedge\Phi_2\wedge\Phi_1^2 + 4\beta_1^2\eta_2\wedge\eta_1\wedge\Phi_1\wedge\Phi_2^2+ 4\beta_2^2\eta_1\wedge\eta_2\wedge\Phi_2^3)$  \\
\hline
5&$\beta_1$-Kenmotsu & $\alpha_2$-Sasakian &0  & 1  & 1 & 0& $3(2\beta_1\alpha_2\Phi_2\wedge\Phi_1^3+ 4\beta_1\alpha_2\Phi_2^2\wedge\Phi_1^2- 24\beta_1\alpha_2\Phi_2\wedge\Phi_1^2\wedge\eta_1\wedge\eta_2+ 4\beta_1^2\eta_1\wedge\eta_2\wedge\Phi_1^3-4\alpha_2^2\Phi_2^3\wedge\Phi_1- 2\alpha_2^2\Phi_2^2\wedge\Phi_1^2+ 16\alpha_2^2\eta_1\wedge\eta_2\wedge\Phi_2^2\wedge\Phi_1+ 2\beta_1\alpha_2\Phi_2^3\wedge\Phi_1- 24\beta_1\alpha_2\Phi_2^2\wedge\Phi_1\wedge\eta_1\wedge\eta_2+ 4\beta_1^2\eta_2\wedge\eta_1\wedge\Phi_1\wedge\Phi_2^2)$  \\
\hline
6&$\beta_1$-Kenmotsu & Cosymplectic & 0 & 0 & 1 & 0 & $3(4\beta_1^2\eta_1\wedge\eta_2\wedge\Phi_1^3+ 4\beta_1^2\eta_2\wedge\eta_1\wedge\Phi_1\wedge\Phi_2^2) $\\
\hline
7&Cosymplectic & $\alpha_2$-Sasakian  & 0 & 1 & 0 & 0 & $3(-4\alpha_2^2\Phi_2^3\wedge\Phi_1+ 16\alpha_2^2\eta_1\wedge\eta_2\wedge\Phi_2^2\wedge\Phi_1+ 16\alpha_2^2\Phi_2^3\wedge\eta_1\wedge\eta_2)$ \\
\hline
8&Cosymplectic &  $\beta_2$-Kenmotsu  & 0 & 0 & 0 & 1 &$ 12\beta_2^2\eta_1\wedge\eta_2\wedge\Phi_2^3$ \\
 \hline
 9&Cosymplectic &  Cosymplectic  & 0 & 0 & 0 & 0 & $0$ \\
 \hline
\end{tabular}

\begin{rem-new}
    1. In Table $4$, cases $3$ and $9$ are the astheno-K\"ahler manifold.\\
    2. A specific instance of Case 3 can be observed in Theorem 4.2 in the paper by Matsuo and Takahashi \cite{km}.
\end{rem-new}

(ii) $\dim_{\mathbbm{R}} M_1 = 5$ and $\dim_{\mathbbm{R}} M_2 = 5$

\vspace{.5cm}
\noindent
\begin{tabular}{|p{.8 cm}|p{2.5cm}|p{2.5cm}|p{.2 cm}|p{.2 cm}|p{.2 cm}|p{.2 cm}|p{6.8 cm}|  }
 \hline
 \multicolumn{8}{|c|}{Table 5: Particular cases} \\
 \hline
 S.No. &$M_1$ & $M_2$& $\alpha_1$ & $\alpha_2$ & $\beta_1$&$\beta_2$& $dd^c \Omega^3 $ \\
 \hline
1&$\alpha_1-$ Sasakian & $\alpha_2$-Sasakian  & 1 & 1 &0  & 0 & $3(-2\alpha_2^2\Phi_2^2\wedge\Phi_1^2 + 16\alpha_2^2\eta_1\wedge\eta_2\wedge\Phi_2^2\wedge\Phi_1 +2\alpha_1^2\Phi_1^2\wedge\Phi_2^2)$\\
\hline
2&$\alpha_1$-Sasakian & $\beta_2$-Kenmotsu  & 1 &0  & 0 &1  &$3(-4\beta_2\alpha_1\Phi_1^2\Phi_2^2 + 8\beta_2\alpha_1\Phi_1^2\wedge\Phi_2\wedge\eta_2\wedge\eta_1 - 4\beta_2^2\eta_1\wedge\eta_2\wedge\Phi_2\wedge\Phi_1^2 -2\beta_2\alpha_1\Phi_1\wedge\Phi_2^3 + 8\beta_2\alpha_1\Phi_1\wedge\Phi_2^2\wedge\eta_2\wedge\eta_1 + 4\beta_2^2\eta_1\wedge\eta_2\wedge\Phi_2^3 +2\alpha_1^2\Phi_1^2\wedge\Phi_2^2)$ \\
\hline
3&$\alpha_1$-Sasakian & Cosymplectic   &1  & 0 & 0 &0  & $6\alpha_1^2\Phi_1^2\wedge\Phi_2^2$ \\
\hline
4&$\beta_1$-Kenmotsu & $\beta_2$-Kenmotsu  & 0 & 0 & 1 &1 & $3(4\beta_1^2\eta_1\wedge\eta_2\wedge\Phi_1^3 -4\beta_2^2\eta_1\wedge\eta_2\wedge\Phi_2\wedge\Phi_1^2 +4\beta_1^2\eta_2\wedge\eta_1\wedge\Phi_1\wedge\Phi_2^2 + 4\beta_2^2\eta_1\wedge\eta_2\wedge\Phi_2^3)$\\
\hline
5&$\beta_1$-Kenmotsu & $\alpha_2$-Sasakian &0  & 1  & 1 & 0&$3(2\beta_1\alpha_2\Phi_2\wedge\Phi_1^3 + 4\beta_1\alpha_2\Phi_2^2\wedge\Phi_1^2 -24\beta_1\alpha_2\Phi_2\wedge\Phi_1^2\wedge\eta_1\wedge\eta_2 + 4\beta_1^2\eta_1\wedge\eta_2\wedge\Phi_1^3 -2\alpha_2^2\Phi_2^2\wedge\Phi_1^2 +16\alpha_2^2\eta_1\wedge\eta_2\wedge\Phi_2^2\wedge\Phi_1 -24\beta_1\alpha_2\Phi_2^2\wedge\Phi_1\wedge\eta_1\wedge\eta_2 -4\beta_1^2\eta_2\wedge\eta_1\wedge\Phi_1\wedge\Phi_2^2)$    \\
\hline
6&$\beta_1$-Kenmotsu & Cosymplectic & 0 & 0 & 1 & 0 &  $3(4\beta_1^2\eta_1\wedge\eta_2\wedge\Phi_1^3 +4\beta_1^2\eta_2\wedge\eta_1\wedge\Phi_1\wedge\Phi_2^2 )$\\
\hline
7&Cosymplectic & $\alpha_2$-Sasakian  & 0 & 1 & 0 & 0 &  $3(-2\alpha_2^2\Phi_2^2\wedge\Phi_1^2 + 16\alpha_2^2\eta_1\wedge\eta_2\wedge\Phi_2^2\wedge\Phi_1) $\\
\hline
8&Cosymplectic &  $\beta_2$-Kenmotsu  & 0 & 0 & 0 & 1 & $ 3(-4\beta_2^2\eta_1\wedge\eta_2\wedge\Phi_2\wedge\Phi_1^2  + 4\beta_2^2\eta_1\wedge\eta_2\wedge\Phi_2^3)$\\
 \hline
 9&Cosymplectic &  Cosymplectic  & 0 & 0 & 0 & 0 & $0$ \\
 \hline
\end{tabular}
\begin{rem-new}
   In Table 5, Case 9 is the astheno-K\"ahler manifold.
\end{rem-new}

 (iii) $\dim_{\mathbbm{R}} M_1 = 7$ and $\dim_{\mathbbm{R}} M_2 = 3$

\vspace{.5cm}
\noindent
\begin{tabular}{|p{.8 cm}|p{2.5cm}|p{2.5cm}|p{.2 cm}|p{.2 cm}|p{.2 cm}|p{.2 cm}|p{6.8 cm}|  }
 \hline
 \multicolumn{8}{|c|}{Table 6: Particular cases} \\
 \hline
 S.No. &$M_1$ & $M_2$& $\alpha_1$ & $\alpha_2$ & $\beta_1$&$\beta_2$& $dd^c \Omega^3 $ \\
 \hline
1&$\alpha_1-$ Sasakian & $\alpha_2$-Sasakian  & 1 & 1 &0  & 0 & $12\alpha_1^2\Phi_1^3\wedge\Phi_2$\\
\hline
2&$\alpha_1$-Sasakian & $\beta_2$-Kenmotsu  & 1 &0  & 0 &1  & $3(-2\beta_2\alpha_1\Phi_1^3\wedge\Phi_2 -4\beta_2\alpha_1\Phi_1^2\wedge\Phi_2^2 + 8\beta_2\alpha_1\Phi_1^2\wedge\Phi_2\wedge\eta_2\wedge\eta_1 -4\beta_2^2\eta_1\wedge\eta_2\wedge\Phi_2\wedge\Phi_1^2 +4\alpha_1^2\Phi_1^3\wedge\Phi_2 -2\beta_2\alpha_1\Phi_1\wedge\Phi_2^3 +8 \beta_2\alpha_1\Phi_1\wedge\Phi_2^2\wedge\eta_2\wedge\eta_1 +4\beta_2^2\eta_1\wedge\eta_2\wedge\Phi_2^3 + 2\alpha_1^2\Phi_1^2\Phi_1^2\wedge\Phi_2^2)$\\
\hline
3&$\alpha_1$-Sasakian & Cosymplectic   &1  & 0 & 0 &0  & $12\alpha_1^2\Phi_1^3\wedge\Phi_2$ \\
\hline
4&$\beta_1$-Kenmotsu & $\beta_2$-Kenmotsu  & 0 & 0 & 1 &1 & $3(4\beta_1^2\eta_1\wedge\eta_2\wedge\Phi_1^3 -4\beta_2^2\eta_1\wedge\eta_2\wedge\Phi_2\wedge\Phi_1^2 + 4\beta_2^2\eta_1\wedge\eta_2\wedge\Phi_2^3 + 4\beta_1^2\eta_2\wedge\eta_1\wedge\Phi_1\wedge\Phi_2^2)$ \\
\hline
5&$\beta_1$-Kenmotsu & $\alpha_2$-Sasakian &0  & 1  & 1 & 0&  $6\beta_1\alpha_2\Phi_2\wedge\Phi_1^3$  \\
\hline
6&$\beta_1$-Kenmotsu & Cosymplectic & 0 & 0 & 1 & 0 &  $12\beta_1^2\eta_1\wedge\eta_2\wedge\Phi_1^3$\\
\hline
7&Cosymplectic & $\alpha_2$-Sasakian  & 0 & 1 & 0 & 0 &  $0$\\
\hline
8&Cosymplectic &  $\beta_2$-Kenmotsu  & 0 & 0 & 0 & 1 &$3(-4\beta_2^2\eta_1\wedge\eta_2\wedge\Phi_2\wedge\Phi_1^2 + 4\beta_2^2\eta_1\wedge\eta_2\wedge\Phi_2^3)$ \\
 \hline
 9&Cosymplectic &  Cosymplectic  & 0 & 0 & 0 & 0 & $0$ \\
 \hline
\end{tabular}
\begin{rem-new}
   1. In Table 6, Cases 7 and 9 are the astheno-K\"ahler manifold.\\
   2. A specific instance of Case 7 can be observed in Theorem 4.2 in the paper by Matsuo and Takahashi \cite{km}.
\end{rem-new}

\noindent
\textbf{Case 4}: For $\dim_{\mathbbm{C}} M = m = 6$, that is, $\dim_{\mathbbm{R}} M = 12$. We have the following possiblities :
\begin{itemize}
          \item [(i)] $\dim_{\mathbbm{R}} M_1 = 3$ \& $\dim_{\mathbbm{R}} M_2 = 9$.
         \item [(ii)]$\dim_{\mathbbm{R}} M_1 = 5$ \& $\dim_{\mathbbm{R}} M_2 = 7$
        \item [(iii)] $\dim_{\mathbbm{R}} M_1 = 7$ \& $\dim_{\mathbbm{R}} M_2 = 5$
       \item [(iv)] $\dim_{\mathbbm{R}} M_1 = 9$ \& $\dim_{\mathbbm{R}} M_2 = 3$.
       \end{itemize}
For astheno-K\"ahler manifold with $m = 6$ , we need to show $$dd^c \Omega^4 =4(dd^c \Omega \wedge \Omega + 3d \Omega \wedge d^c \Omega) \wedge \Omega^2 = 0$$
By (\ref{p6}) and (\ref{p7}),we have
\begin{eqnarray*}
 dd^c \Omega^4 &=&4\Big[2(\beta_1 \alpha_2 \Phi_2 \wedge \Phi_1^2+\beta_1 \alpha_2 \Phi_2^2 \wedge\Phi_1-2 \beta_1 \alpha_2 \Phi_2 \wedge \Phi_1 \wedge \eta_1 \wedge \eta_2 +2 \beta_1^2 \eta_2 \wedge \eta_1 \wedge \Phi_1^2 \notag \\
&&+2 \beta_1^2 \eta_2 \wedge \eta_1 \wedge \Phi_1 \wedge \Phi_2 -4 \beta_1^2 \eta_2 \wedge \eta_1\wedge \Phi_1\wedge \eta_1 \wedge \eta_2  -\beta_2 \alpha_1 \Phi_1^2 \wedge \Phi_2-\beta_2 \alpha_1 \Phi_1 \wedge \Phi_2^2 \notag\\
&&+2 \beta_2 \alpha_1 \Phi_1 \wedge \Phi_2 \wedge \eta_1 \wedge \eta_2 -2 \beta_2^2 \eta_1 \wedge \eta_2 \wedge \Phi_2 \wedge \Phi_1-2 \beta_2^2 \wedge \eta_1 \wedge \eta_2 \wedge \Phi_2^{2}\notag  \\
&& +4 \beta_2^2 \eta_1 \wedge \eta_2 \wedge \Phi_2 \wedge \eta_1 \wedge \eta_2 +\alpha_1^2 \Phi_1^3 +\alpha_1^2 \Phi_1^2 \wedge \Phi_2-2 \alpha_1^2 \Phi_1^2 \wedge \eta_1 \wedge \eta_2 -\alpha_2^2 \Phi_2^2 \wedge \Phi_1\notag \\
&& -\alpha_2^2 \Phi_2^3+2 \alpha_2^2 \Phi_2^2 \wedge \eta_1 \wedge \eta_2) + 3\times4(\beta_1^2 \eta_1 \wedge \eta_2 \wedge \Phi_1^2-\beta_1 \alpha_2 \eta_1 \wedge \eta_2 \wedge \Phi_1 \wedge \Phi_2 \notag\\
&&-\beta_2^2 \eta_2 \wedge \eta_1 \wedge \Phi_2^2 +\beta_2 \alpha_1 \eta_2 \wedge \eta_1\wedge\Phi_2 \wedge \Phi_1  +\alpha_1 \beta_2 \eta_2 \wedge \eta_1 \wedge \Phi_1 \wedge \Phi_2-\alpha_1^2 \eta_2 \wedge \eta_1 \wedge \Phi_1^2 \notag \\
&&-\beta_1 \alpha_2 \eta_1 \wedge \eta_2 \wedge \Phi_2 \wedge \Phi_1+\alpha_2^2 \eta_1\wedge \eta_2 \wedge \Phi_2^2) \Big]\wedge\Big[\Phi_1+\Phi_2-2 \eta_1 \wedge \eta_2\Big]^2 
\end{eqnarray*}

\begin{eqnarray*}
     &=&4\Big[2\beta_1\alpha_2\Phi_2\wedge\Phi_1^4 + 6\beta_1\alpha_2\Phi_2^2\wedge\Phi_1^3 -36\beta_1\alpha_2\Phi_2\wedge\Phi_1^3\wedge\eta_1\wedge\eta_2 + 8\beta_1^2\eta_1\wedge\eta_2\wedge\Phi_1^4 + \\
    && 12\beta_1^2\eta_1\wedge\eta_2\wedge\Phi_1^3\wedge\Phi_2 -2\beta_2\alpha_1\Phi_1^4\wedge\Phi_2 -6\beta_2\alpha_1\Phi_1^3\wedge\Phi_2^2 -4\beta_2^2\eta_1\wedge\eta_2\wedge\Phi_2\wedge\Phi_1^3 +2\alpha_1^2\Phi_1^5 \\
    && +6\alpha_1^2\Phi_1^4\wedge\Phi_2 - 2\alpha_2^2\Phi_2^2\wedge\Phi_1^3 -6\alpha_2^2\Phi_2^3\wedge\Phi_1^2 + 24\alpha_2^2\Phi_2^2\wedge\Phi_1^2\wedge\eta_1\wedge\eta_2 +6\beta_1\alpha_2\Phi_2^3\wedge\Phi_1^2 \\
    &&+ 2\beta_1\alpha_2\Phi_2^4\wedge\Phi_1 -36\beta_1\alpha_2\Phi_2^3\wedge\Phi_1\wedge\eta_1\wedge\eta_2 +4\beta_1^2\eta_2\wedge\eta_1\wedge\Phi_1\wedge\Phi_2^3 -6\beta_2\alpha_1\Phi_1^2\wedge\Phi_2^3 \\
    && -2\beta_2\alpha_1\Phi_1\wedge\Phi_2^4 +12\beta_2^2\eta_1\wedge\eta_2\wedge\Phi_2^3\wedge\Phi_1 + 8\beta_2^2\eta_1\wedge\eta_2\wedge\Phi_2^4 + 6\alpha_1^2\Phi_1^3\wedge\Phi_2^2 + 2\alpha_1^2\Phi_1^2\wedge\Phi_2^3 \\
    &&- 6\alpha_2^2\Phi_2^4\wedge\Phi_1 -2\alpha_2^2\Phi_2^5 + 24\alpha_2^2\Phi_2^4\wedge\eta_1\wedge\eta_2 -60\beta_1\alpha_2\Phi_2^2\wedge\Phi_1^2\wedge\eta_1\wedge\eta_2 \\
    &&+ 8\beta_2\alpha_1\eta_2\wedge\eta_1\wedge\Phi_2^2\wedge\Phi_1^2 
    - 16\alpha_1^2\eta_1\wedge\eta_2\wedge\Phi_1^3\wedge\Phi_2 + 48\alpha_2^2\Phi_2^3\wedge\Phi_1\wedge\eta_1\wedge\eta_2\Big].
\end{eqnarray*}

 (i) $\dim_{\mathbbm{R}} M_1 = 3$ \& $\dim_{\mathbbm{R}} M_2 = 9$\\
 \vspace{.5cm}
\noindent
\begin{tabular}{|p{.8 cm}|p{2.5cm}|p{2.5cm}|p{.2 cm}|p{.2 cm}|p{.2 cm}|p{.2 cm}|p{6.8 cm}|  }
 \hline
 \multicolumn{8}{|c|}{Table 7: Particular cases} \\
 \hline
 S.No. &$M_1$ & $M_2$& $\alpha_1$ & $\alpha_2$ & $\beta_1$&$\beta_2$& $dd^c \Omega^4 $ \\
 \hline
1&$\alpha_1$-Sasakian & $\alpha_2$-Sasakian  & 1 & 1 &0  & 0 & $4(-6\alpha_2^2\Phi_2^4\wedge\Phi_1 + 24\alpha_2^2\Phi_2^4\wedge\eta_1\wedge\eta_2 + 48\alpha_2^2\Phi_2^3\wedge\Phi_1\wedge\eta_1\wedge\eta_2)$\\
\hline
2&$\alpha_1$-Sasakian & $\beta_2$-Kenmotsu  & 1 &0  & 0 &1  &$4(-2\beta_2\alpha_1\Phi_1\wedge\Phi_2^4 +12\beta_2^2\eta_1\wedge\eta_2\wedge\Phi_2^3\wedge\Phi_1 +8\beta_2^2\eta_1\wedge\eta_2\wedge\Phi_2^4)$ \\
\hline
3&$\alpha_1$-Sasakian & Cosymplectic   &1  & 0 & 0 &0  & $0$ \\
\hline
4&$\beta_1$-Kenmotsu & $\beta_2$-Kenmotsu  & 0 & 0 & 1 &1 & $4(8\beta_1^2\eta_1\wedge\eta_2\wedge\Phi_1^4 +12\beta_1^2\eta_1\wedge\eta_2\wedge\Phi_1^3\wedge\Phi_2 -4\beta_2^2\eta_1\wedge\eta_2\wedge\Phi_2\wedge\Phi_1^3 +4\beta_1^2\eta_2\wedge\eta_1\wedge\Phi_1\wedge\Phi_2^3 +12\beta_2^2\eta_1\wedge\eta_2\wedge\Phi_2^3\wedge\Phi_1 +8\beta_2^2\eta_1\wedge\eta_2\wedge\Phi_2^4)$\\
\hline
5&$\beta_1$-Kenmotsu & $\alpha_2$-Sasakian &0  & 1  & 1 & 0&  $4(2\beta_1\alpha_2\Phi_2\wedge\Phi_1^4 + 6\beta_1\alpha_2\Phi_2^2\wedge\Phi_1^3 -36\beta_1\alpha_2\Phi_2\wedge\Phi_1^3\wedge\eta_1\wedge\eta_2 + 8\beta_1^2\eta_1\wedge\eta_2\wedge\Phi_1^4+ 12\beta_1^2\eta_1\wedge\eta_2\wedge\Phi_1^4 - 2\alpha_2^2\Phi_2^2\wedge\Phi_1^3 -6\alpha_2^2\Phi_2^3\wedge\Phi_1^2 + 24\alpha_2^2\Phi_2^2\wedge\Phi_1^2\wedge\eta_1\wedge\eta_2 +6\beta_1\alpha_2\Phi_2^3\wedge\Phi_1^2 + 2\beta_1\alpha_2\Phi_2^4\wedge\Phi_1 -36\beta_1\alpha_2\Phi_2^3\wedge\Phi_1\wedge\eta_1\wedge\eta_2 +4\beta_1^2\eta_2\wedge\eta_1\wedge\Phi_1\wedge\Phi_2^3 - 6\alpha_2^2\Phi_2^4\wedge\Phi_1 + 24\alpha_2^2\Phi_2^4\wedge\eta_1\wedge\eta_2 -60\beta_1\alpha_2\Phi_2^2\wedge\Phi_1^2\wedge\eta_1\wedge\eta_2+ 48\alpha_2^2\Phi_2^3\wedge\Phi_1\wedge\eta_1\wedge\eta_2)$  \\
\hline
6&$\beta_1$-Kenmotsu & Cosymplectic & 0 & 0 & 1 & 0 & $4(8\beta_1^2\eta_1\wedge\eta_2\wedge\Phi_1^4 +12\beta_1^2\eta_1\wedge\eta_2\wedge\Phi_1^3\wedge\Phi_2 +4\beta_1^2\eta_2\wedge\eta_1\wedge\Phi_1\wedge\Phi_2^3 )$ \\
\hline
7&Cosymplectic & $\alpha_2$-Sasakian  & 0 & 1 & 0 & 0 & $4(-6\alpha_2^2\Phi_2^4\wedge\Phi_1 + 24\alpha_2^2\Phi_2^4\wedge\eta_1\wedge\eta_2 + 48\alpha_2^2\Phi_2^3\wedge\Phi_1\wedge\eta_1\wedge\eta_2)$ \\
\hline
8&Cosymplectic &  $\beta_2$-Kenmotsu  & 0 & 0 & 0 & 1 &$4(12\beta_2^2\eta_1\wedge\eta_2\wedge\Phi_2^3\wedge\Phi_1 +8\beta_2^2\eta_1\wedge\eta_2\wedge\Phi_2^4)$ \\
 \hline
 9&Cosymplectic &  Cosymplectic  & 0 & 0 & 0 & 0 &  $0$\\
 \hline
\end{tabular}
\begin{rem-new}
   1. In Table 7, Cases 3 and 9 are the astheno-K\"ahler manifold.\\
   2. A specific instance of Case 3 can be observed in Theorem 4.2 in the paper by Matsuo and Takahashi \cite{km}.
\end{rem-new}

(ii)$\dim_{\mathbbm{R}} M_1 = 5$ \& $\dim_{\mathbbm{R}} M_2 = 7$\\
\vspace{.5cm}
\noindent
\begin{tabular}{|p{.8 cm}|p{2.5cm}|p{2.5cm}|p{.2 cm}|p{.2 cm}|p{.2 cm}|p{.2 cm}|p{6.8 cm}|  }
 \hline
 \multicolumn{8}{|c|}{Table 8: Particular cases} \\
 \hline
 S.No. &$M_1$ & $M_2$& $\alpha_1$ & $\alpha_2$ & $\beta_1$&$\beta_2$& $dd^c \Omega^4 $ \\
 \hline
1&$\alpha_1$-Sasakian & $\alpha_2$-Sasakian  & 1 & 1 &0  & 0 & $4(2\alpha_1^2\Phi_1^2\wedge\Phi_2^3- 6\alpha_2^2\Phi_2^3\wedge\Phi_1^2 + 24\alpha_2^2\Phi_2^2\wedge\Phi_1^2\wedge\eta_1\wedge\eta_2 + 48\alpha_2^2\Phi_2^3\wedge\Phi_1\wedge\eta_1\wedge\eta_2 )$\\
\hline
2&$\alpha_1$-Sasakian & $\beta_2$-Kenmotsu  & 1 &0  & 0 &1  &$4(-6\beta_2\alpha_1\Phi_1^2\wedge\Phi_2^3 -2\beta_2\alpha_1\Phi_1\wedge\Phi_2^4 +12\beta_2^2\eta_1\wedge\eta_2\wedge\Phi_2^3\wedge\Phi_1 + 8\beta_2^2\eta_1\wedge\eta_2\wedge\Phi_2^4 + 2\alpha_1^2\Phi_1^2\wedge\Phi_2^3+ 8\beta_2\alpha_1\eta_2\wedge\eta_1\wedge\Phi_2^2\wedge\Phi_1^2 )$ \\
\hline
3&$\alpha_1$-Sasakian & Cosymplectic   &1  & 0 & 0 &0  &$8\alpha_1^2\Phi_1^2\wedge\Phi_2^3$ \\
\hline
4&$\beta_1$-Kenmotsu & $\beta_2$-Kenmotsu  & 0 & 0 & 1 &1 & $4(8\beta_1^2\eta_1\wedge\eta_2\wedge\Phi_1^4 +12\beta_1^2\eta_1\wedge\eta_2\wedge\Phi_1^3\wedge\Phi_2 -4\beta_2^2\eta_1\wedge\eta_2\wedge\Phi_2\wedge\Phi_1^3 +4\beta_1^2\eta_2\wedge\eta_1\wedge\Phi_1\wedge\Phi_2^3 +12\beta_2^2\eta_1\wedge\eta_2\wedge\Phi_2^3\wedge\Phi_1 +8\beta_2^2\eta_1\wedge\eta_2\wedge\Phi_2^4)$\\
\hline
5&$\beta_1$-Kenmotsu & $\alpha_2$-Sasakian &0  & 1  & 1 & 0& $4(2\beta_1\alpha_2\Phi_2\wedge\Phi_1^4 + 6\beta_1\alpha_2\Phi_2^2\wedge\Phi_1^3 -36\beta_1\alpha_2\Phi_2\wedge\Phi_1^3\wedge\eta_1\wedge\eta_2 + 8\beta_1^2\eta_1\wedge\eta_2\wedge\Phi_1^4 + 12\beta_1^2\eta_1\wedge\eta_2\wedge\Phi_1^4- 2\alpha_2^2\Phi_2^2\wedge\Phi_1^3 -6\alpha_2^2\Phi_2^3\wedge\Phi_1^2 + 24\alpha_2^2\Phi_2^2\wedge\Phi_1^2\wedge\eta_1\wedge\eta_2 +6\beta_1\alpha_2\Phi_2^3\wedge\Phi_1^2-36\beta_1\alpha_2\Phi_2^3\wedge\Phi_1\wedge\eta_1\wedge\eta_2 +4\beta_1^2\eta_2\wedge\eta_1\wedge\Phi_1\wedge\Phi_2^3 -60\beta_1\alpha_2\Phi_2^2\wedge\Phi_1^2\wedge\eta_1\wedge\eta_2+ 48\alpha_2^2\Phi_2^3\wedge\Phi_1\wedge\eta_1\wedge\eta_2 )$   \\
\hline
6&$\beta_1$-Kenmotsu & Cosymplectic & 0 & 0 & 1 & 0 &  $4(8\beta_1^2\eta_1\wedge\eta_2\wedge\Phi_1^4 +12\beta_1^2\eta_1\wedge\eta_2\wedge\Phi_1^3\wedge\Phi_2 +4\beta_1^2\eta_2\wedge\eta_1\wedge\Phi_1\wedge\Phi_2^3 )$ \\
\hline
7&Cosymplectic & $\alpha_2$-Sasakian  & 0 & 1 & 0 & 0 & $4(- 6\alpha_2^2\Phi_2^3\wedge\Phi_1^2 + 24\alpha_2^2\Phi_2^2\wedge\Phi_1^2\wedge\eta_1\wedge\eta_2 + 48\alpha_2^2\Phi_2^3\wedge\Phi_1\wedge\eta_1\wedge\eta_2)$ \\
\hline
8&Cosymplectic &  $\beta_2$-Kenmotsu  & 0 & 0 & 0 & 1 & $4(12\beta_2^2\eta_1\wedge\eta_2\wedge\Phi_2^3\wedge\Phi_1 +8\beta_2^2\eta_1\wedge\eta_2\wedge\Phi_2^4)$ \\
 \hline
 9&Cosymplectic &  Cosymplectic  & 0 & 0 & 0 & 0 & $0$ \\
 \hline
\end{tabular}

\vspace{-.5cm}
(iii) $\dim_{\mathbbm{R}} M_1 = 7$ \& $\dim_{\mathbbm{R}} M_2 = 5$\\
\vspace{.5cm}
\noindent
\begin{tabular}{|p{.8 cm}|p{2.5cm}|p{2.5cm}|p{.2 cm}|p{.2 cm}|p{.2 cm}|p{.2 cm}|p{6.8 cm}|  }
 \hline
 \multicolumn{8}{|c|}{Table 9: Particular cases} \\
 \hline
 S.No. &$M_1$ & $M_2$& $\alpha_1$ & $\alpha_2$ & $\beta_1$&$\beta_2$& $dd^c \Omega^4 $ \\
 \hline
1&$\alpha_1$-Sasakian & $\alpha_2$-Sasakian  & 1 & 1 &0  & 0 & $4(- 2\alpha_2^2\Phi_2^2\wedge\Phi_1^3 + 24\alpha_2^2\Phi_2^2\wedge\Phi_1^2\wedge\eta_1\wedge\eta_2  + 6\alpha_1^2\Phi_1^3\wedge\Phi_2^2 - 16\alpha_1^2\eta_1\wedge\eta_2\wedge\Phi_1^3\wedge\Phi_2 )$\\
\hline
2&$\alpha_1$-Sasakian & $\beta_2$-Kenmotsu  & 1 &0  & 0 &1  & $4(-6\beta_2\alpha_1\Phi_1^3\wedge\Phi_2^2 -4\beta_2^2\eta_1\wedge\eta_2\wedge\Phi_2\wedge\Phi_1^3-6\beta_2\alpha_1\Phi_1^2\wedge\Phi_2^3 -2\beta_2\alpha_1\Phi_1\wedge\Phi_2^4 +12\beta_2^2\eta_1\wedge\eta_2\wedge\Phi_2^3\wedge\Phi_1 + 8\beta_2^2\eta_1\wedge\eta_2\wedge\Phi_2^4 + 6\alpha_1^2\Phi_1^3\wedge\Phi_2^2 + 2\alpha_1^2\Phi_1^2\wedge\Phi_2^3 + 8\beta_2\alpha_1\eta_2\wedge\eta_1\wedge\Phi_2^2\wedge\Phi_1^2 
    - 16\alpha_1^2\eta_1\wedge\eta_2\wedge\Phi_1^3\wedge\Phi_2)$\\
\hline
3&$\alpha_1$-Sasakian & Cosymplectic   &1  & 0 & 0 &0  & $4(6\alpha_1^2\Phi_1^3\wedge\Phi_2^2 - 16\alpha_1^2\eta_1\wedge\eta_2\wedge\Phi_1^3\wedge\Phi_2)$ \\
\hline
4&$\beta_1$-Kenmotsu & $\beta_2$-Kenmotsu  & 0 & 0 & 1 &1 & $4(8\beta_1^2\eta_1\wedge\eta_2\wedge\Phi_1^4 +12\beta_1^2\eta_1\wedge\eta_2\wedge\Phi_1^3\wedge\Phi_2 -4\beta_2^2\eta_1\wedge\eta_2\wedge\Phi_2\wedge\Phi_1^3 +4\beta_1^2\eta_2\wedge\eta_1\wedge\Phi_1\wedge\Phi_2^3 +12\beta_2^2\eta_1\wedge\eta_2\wedge\Phi_2^3\wedge\Phi_1 +8\beta_2^2\eta_1\wedge\eta_2\wedge\Phi_2^4)$\\
\hline
5&$\beta_1$-Kenmotsu & $\alpha_2$-Sasakian &0  & 1  & 1 & 0&  $4(2\beta_1\alpha_2\Phi_2\wedge\Phi_1^4 + 6\beta_1\alpha_2\Phi_2^2\wedge\Phi_1^3 -36\beta_1\alpha_2\Phi_2\wedge\Phi_1^3\wedge\eta_1\wedge\eta_2 + 8\beta_1^2\eta_1\wedge\eta_2\wedge\Phi_1^4 + 12\beta_1^2\eta_1\wedge\eta_2\wedge\Phi_1^3\wedge\Phi_2- 2\alpha_2^2\Phi_2^2\wedge\Phi_1^3 + 24\alpha_2^2\Phi_2^2\wedge\Phi_1^2\wedge\eta_1\wedge\eta_2-60\beta_1\alpha_2\Phi_2^2\wedge\Phi_1^2\wedge\eta_1\wedge\eta_2  )$  \\
\hline
6&$\beta_1$-Kenmotsu & Cosymplectic & 0 & 0 & 1 & 0 &  $4(8\beta_1^2\eta_1\wedge\eta_2\wedge\Phi_1^4 +12\beta_1^2\eta_1\wedge\eta_2\wedge\Phi_1^3\wedge\Phi_2 )$\\
\hline
7&Cosymplectic & $\alpha_2$-Sasakian  & 0 & 1 & 0 & 0 & $4(- 2\alpha_2^2\Phi_2^2\wedge\Phi_1^3 + 24\alpha_2^2\Phi_2^2\wedge\Phi_1^2\wedge\eta_1\wedge\eta_2 )$ \\
\hline
8&Cosymplectic &  $\beta_2$-Kenmotsu  & 0 & 0 & 0 & 1 &$4(-4\beta_2^2\eta_1\wedge\eta_2\wedge\Phi_2\wedge\Phi_1^3 +12\beta_2^2\eta_1\wedge\eta_2\wedge\Phi_2^3\wedge\Phi_1 +8\beta_2^2\eta_1\wedge\eta_2\wedge\Phi_2^4)$ \\
 \hline
 9&Cosymplectic &  Cosymplectic  & 0 & 0 & 0 & 0 & $0$ \\
 \hline
\end{tabular}
\begin{rem-new}
   In Table 8, Case 9 is the astheno-K\"ahler manifold.
\end{rem-new}

\begin{rem-new}
   In Table 9, Case 9 is the astheno-K\"ahler manifold.
\end{rem-new}
(iv) $\dim_{\mathbbm{R}} M_1 = 9$ \& $\dim_{\mathbbm{R}} M_2 = 3$\\
\begin{tabular}{|p{.8 cm}|p{2.5cm}|p{2.5cm}|p{.2 cm}|p{.2 cm}|p{.2 cm}|p{.2 cm}|p{6.8 cm}|  }
 \hline
 \multicolumn{8}{|c|}{Table 10: Particular cases} \\
 \hline
 S.No. &$M_1$ & $M_2$& $\alpha_1$ & $\alpha_2$ & $\beta_1$&$\beta_2$& $dd^c \Omega^4 $ \\
 \hline
1&$\alpha_1$-Sasakian & $\alpha_2$-Sasakian  & 1 & 1 &0  & 0 & $4(6\alpha_1^2\Phi_1^4\wedge\Phi_2- 16\alpha_1^2\eta_1\wedge\eta_2\wedge\Phi_1^3\wedge\Phi_2)$\\
\hline
2&$\alpha_1$-Sasakian & $\beta_2$-Kenmotsu  & 1 &0  & 0 &1  &$4(-2\beta_2\alpha_1\Phi_1^4\wedge\Phi_2 -6\beta_2\alpha_1\Phi_1^3\wedge\Phi_2^2 -4\beta_2^2\eta_1\wedge\eta_2\wedge\Phi_2\wedge\Phi_1^3 +6\alpha_1^2\Phi_1^4\wedge\Phi_2-6\beta_2\alpha_1\Phi_1^2\wedge\Phi_2^3  -2\beta_2\alpha_1\Phi_1\wedge\Phi_2^4 +12\beta_2^2\eta_1\wedge\eta_2\wedge\Phi_2^3\wedge\Phi_1 + 8\beta_2^2\eta_1\wedge\eta_2\wedge\Phi_2^4 + 6\alpha_1^2\Phi_1^3\wedge\Phi_2^2 + 2\alpha_1^2\Phi_1^2\wedge\Phi_2^3+ 8\beta_2\alpha_1\eta_2\wedge\eta_1\wedge\Phi_2^2\wedge\Phi_1^2 
    - 16\alpha_1^2\eta_1\wedge\eta_2\wedge\Phi_1^3\wedge\Phi_2)$ \\
\hline
3&$\alpha_1$-Sasakian & Cosymplectic   &1  & 0 & 0 &0  & $4(6\alpha_1^2\Phi_1^4\wedge\Phi_2- 16\alpha_1^2\eta_1\wedge\eta_2\wedge\Phi_1^3\wedge\Phi_2)$ \\
\hline
4&$\beta_1$-Kenmotsu & $\beta_2$-Kenmotsu  & 0 & 0 & 1 &1 & $4(8\beta_1^2\eta_1\wedge\eta_2\wedge\Phi_1^4 +12\beta_1^2\eta_1\wedge\eta_2\wedge\Phi_1^3\wedge\Phi_2 -4\beta_2^2\eta_1\wedge\eta_2\wedge\Phi_2\wedge\Phi_1^3 +4\beta_1^2\eta_2\wedge\eta_1\wedge\Phi_1\wedge\Phi_2^3 +12\beta_2^2\eta_1\wedge\eta_2\wedge\Phi_2^3\wedge\Phi_1 +8\beta_2^2\eta_1\wedge\eta_2\wedge\Phi_2^4)$\\
\hline
5&$\beta_1$-Kenmotsu & $\alpha_2$-Sasakian &0  & 1  & 1 & 0& $4(2\beta_1\alpha_2\Phi_2\wedge\Phi_1^4-36\beta_1\alpha_2\Phi_2\wedge\Phi_1^3\wedge\eta_1\wedge\eta_2 + 8\beta_1^2\eta_1\wedge\eta_2\wedge\Phi_1^4 + 12\beta_1^2\eta_1\wedge\eta_2\wedge\Phi_1^3\wedge\Phi_2)$   \\
\hline
6&$\beta_1$-Kenmotsu & Cosymplectic & 0 & 0 & 1 & 0 & $4(8\beta_1^2\eta_1\wedge\eta_2\wedge\Phi_1^4 +12\beta_1^2\eta_1\wedge\eta_2\wedge\Phi_1^3\wedge\Phi_2)$ \\
\hline
7&Cosymplectic & $\alpha_2$-Sasakian  & 0 & 1 & 0 & 0 & $0$ \\
\hline
8&Cosymplectic &  $\beta_2$-Kenmotsu  & 0 & 0 & 0 & 1 & $4(-4\beta_2^2\eta_1\wedge\eta_2\wedge\Phi_2\wedge\Phi_1^3 +12\beta_2^2\eta_1\wedge\eta_2\wedge\Phi_2^3\wedge\Phi_1 +8\beta_2^2\eta_1\wedge\eta_2\wedge\Phi_2^4)$ \\
 \hline
 9&Cosymplectic &  Cosymplectic  & 0 & 0 & 0 & 0 & $0$ \\
 \hline
\end{tabular}
\begin{rem-new}
   1. In Table 10, Cases 7 and 9 are the astheno-K\"ahler manifold.\\
   2. A specific instance of Case 7 can be observed in Theorem 4.2 in the paper by Matsuo and Takahashi \cite{km}.
\end{rem-new}

\begin{rem-new}
   From \cite{capursi} the case when the manifold is a product of two cosymplectic manifolds, it will always be a astheno-K\"ahler manifold.
\end{rem-new}
\begin{rem-new} By the above study and using the result of Matsuo and Takahashi {\rm\cite[Theorem 4.2]{km}} and Theorem 2.11, we can say that, if $\dim_{\mathbbm{R}} M_1\geq 5$ and $\dim_{\mathbbm{R}} M_2 \geq 5$ then the only possibility of product manifold $M$ to be astheno-K\"ahler is both $M_1$, $M_2$ should be cosymplectic.  
Also, we can observe that for $\dim_{\mathbbm{C}} M=m\geq 3$, a product of $\alpha$-Sasakian and cosymplectic is astheno-K\"ahler only when the real dimension of $\alpha$-Sasakian is $3$.

\end{rem-new}
\noindent
\textbf{Declaration:}  We declare that no conflicts of interest are associated with this research work and there has been no significant financial support for this work.  We certify that the submission is original work.

\end{document}